\documentclass[12pt]{article}
\usepackage{amsmath,amsthm,amssymb}
\usepackage{multirow}

\newtheorem{theorem}{Theorem}[section]
\newtheorem{lemma}[theorem]{Lemma}
\newtheorem{corollary}[theorem]{Corollary}

\numberwithin{equation}{section}

\newcommand{\ssection}[1]{
     \section{\normalsize\bf #1}}

\begin{document}

\title {\large\sc Solution of the Dirichlet problem by a finite difference
analog of the boundary integral equation
 }

\author{
{\normalsize J. Thomas Beale} \\
{\normalsize{\em Department of Mathematics, Duke University, Box 90320}}\\
{\normalsize{\em Durham, North Carolina 27708, U.S.A.}}\\
\and
{\normalsize Wenjun Ying} \\
{\normalsize{\em School of Mathematical Sciences, MOE-LSC and Institute of Natural Sciences}}\\
{\normalsize{\em Shanghai Jiao Tong University, Minhang, Shanghai 200240, P. R. China}}\\
}

\date{}

\maketitle

{\footnotetext {
This work was supported in part by NSF 
grant DMS-1312654 and NSF-China grant DMS-11771290. } }


\newcommand{\beq}{\begin{equation}}
\newcommand{\eeq}{\end{equation}}
\newcommand\p{\,+\,}
\newcommand\lee{\,\leq\,}
\newcommand\gee{\,\geq\,}
\newcommand\m{\,-\,}
\newcommand\eq{\,=\,}
\newcommand{\eps}{\varepsilon}
\newcommand{\sig}{\sigma}
\newcommand{\pa}{\partial}
\newcommand{\lilhalf}{{\textstyle \frac12}}
\newcommand{\lilth}{{\textstyle \frac32}}
\newcommand{\mm} {_{\max}}
\newcommand{\np}{{n+1}}
\newcommand{\nm}{{n-1}}
\newcommand{\nhalf}{{n+1/2}}
\newcommand{\ex}{{ex}}
\newcommand{\NN}{{\mathcal N}}
\newcommand{\R}{{\mathbb R}}
\newcommand{\Z}{{\mathbb Z}}
\newcommand{\laph}{\Delta_h}
\newcommand{\Ptw}{{\tilde P}}
\newcommand{\fcns}{{\mathcal F}(\Omega_h)}
\newcommand{\phihat}{{\hat \varphi}}
\newcommand{\fhat}{{\hat f}}
\newcommand{\ghat}{{\hat g}}
\newcommand{\ip}{{i+1}}
\newcommand{\im}{{i-1}}
\newcommand{\iph}{{i+1/2}}
\newcommand{\imh}{{i-1/2}}
\newcommand{\jp}{{j+1}}
\newcommand{\jm}{{j-1}}
\newcommand{\jph}{{j+1/2}}
\newcommand{\jmh}{{j-1/2}}

\newcommand{\kp}{{k+1}}
\newcommand{\km}{{k-1}}
\newcommand{\kph}{{k+1/2}}
\newcommand{\kmh}{{k-1/2}}
\newcommand{\om}{\Omega}

\newcommand{\xtw}{\tilde x}
\newcommand{\ytw}{\tilde y}
\newcommand{\utw}{\tilde u}
\newcommand{\vtw}{\tilde v}
\newcommand{\wtw}{\tilde w}
\newcommand{\phitw}{\tilde \varphi}

\newcommand{\bbox}{{\cal B}}
\newcommand{\boxh}{{\cal B}_h}
\newcommand{\omp}{\om^+}
\newcommand{\omm}{\om^-}
\newcommand{\omhp}{\om_h^+}
\newcommand{\omhm}{\om_h^-}
\newcommand{\omhpm}{\om_h^\pm}
\newcommand{\ompcl}{\overline \omp}
\newcommand{\ommcl}{\overline \omm}
\newcommand{\omhpcl}{\overline \omhp}
\newcommand{\omhmcl}{\overline \omhm}
\newcommand{\omhpmcl}{\overline \omhpm}

\newcommand{\fcnspeq}{{\cal F}_{\#}(\omhpcl)}
\newcommand{\fcnseq}{{\cal F}_{\#}(\Gamma_h^1)}
\newcommand{\fcnscut}{{\cal F}(\Gamma_h^0)}
\newcommand{\fcnsgam}{{\cal F}(\Gamma_h^1)}
\newcommand{\fcnsp}{{\cal F}(\omhpcl)}
\newcommand{\fcnsm}{{\cal F}(\omhmcl)}
\newcommand{\fcnspm}{{\cal F}(\omhpmcl)}
\newcommand{\fcnspz}{{\cal F}_0(\omhpcl)}
\newcommand{\ints}{\cal I}
\newcommand{\intsh}{{\cal I}_h}
\newcommand{\intscut}{{\cal I}^c}
\newcommand{\intshcut}{\intscut_h}
\newcommand{\intsp}{{\ints^+}}
\newcommand{\intsm}{{\ints^-}}

\newcommand{\intshp}{{\cal I}^+_h}
\newcommand{\intshm}{{\cal I}^-_h}

\newcommand{\A}{{\cal A}}
\newcommand{\B}{{\cal B}}

\newcommand{\fbar}{\overline f}

\begin{abstract}
Several important problems in partial differential equations can be formulated as integral equations.  Often the integral operator defines the solution of an elliptic problem with 
specified jump conditions at an interface.  In principle the integral equation can be
solved by replacing the integral operator with a finite difference calculation on a regular grid.  A practical method of this type has been developed by the second author.  In this paper we prove the validity of a simplified version of this method for the Dirichlet problem in a general domain in $\R^2$ or $\R^3$.  Given a boundary value, we solve for a discrete version of the density of the double layer potential using a low order interface method.  It produces the Shortley-Weller solution for the unknown harmonic function with accuracy $O(h^2)$.
We prove the unique solvability for the density, with bounds in norms based on the energy or Dirichlet
norm, using techniques which mimic those of exact potentials.  The analysis reveals that this crude method maintains much of the mathematical structure of the classical integral equation.  Examples are included.

\medskip



\end{abstract} 




\ssection{Introduction}

In classical potential theory the Dirichlet problem on a bounded domain $\omp$ in $\R^2$ or $\R^3$
with boundary $\Gamma$,
\beq \Delta u^+ \eq 0 \mbox{\;\;on\;} \omp \,, \qquad u^+  = g  \mbox{\;\;on\;} \Gamma \eeq
is expressed as an integral equation
$  (\lilhalf + K)\varphi \eq g  $,
where $K$ is the integral operator  for the double layer potential and $\varphi$ is an unknown dipole density function on $\Gamma$.  The double layer potential determined by $\varphi$ is
a pair of harmonic functions, $u^+$ on $\omp$ and $u^-$ on the exterior domain
$\omm = \R^d - \ompcl$, such that
$u^+ - u^- = \varphi$ and $\pa u^+/\pa n = \pa u^-/\pa n$ on $\Gamma$, with
$u^+ = (\lilhalf + K)\varphi \equiv A\varphi$ on $\Gamma$ and
$u^- = (-\lilhalf + K)\varphi \equiv B\varphi$ on $\Gamma$, so that $A = B + I$.  Thus the solution of (1.1) is obtained by inverting $A$.  More specifically, it was proved in
\cite{cost,SW} that $\|B\| < 1$ with certain choices of norms;
the invertibility of $A$ follows from this, along with an estimate for $A^{-1}$.
In \cite{SW} it is shown that $B$ is a contraction if we define $\|\varphi\|$ by
\beq \|\varphi\|^2 \eq \int_{\omp} |\nabla v^+|^2 \p \int_{\omm} |\nabla v^-|^2  \eeq
where $\Delta v^+ = 0$ on $\omp$, $\Delta v^- = 0$ on $\omm$, and
$v^\pm = \varphi$ on $\Gamma$.
This point of view gives a natural setting for the integral equation approach and discrete approximations.
In the present work we use this setting to prove the bounded solvability of a finite difference method which mimics the integral equation.

The double layer potential can be thought of as the solution to an interface problem
\beq \Delta u^\pm = 0 \mbox{\;\;on\;} \om^\pm \,, \qquad [u] = \varphi  \mbox{\;\;on\;} \Gamma\,,
        \qquad [\pa u/\pa n] = 0  \mbox{\;\;on\;} \Gamma
\eeq
where $[u] = u^+ - u^-$.  This problem can be solved by a finite difference method such as the immersed interface method \cite{liitobook,mayo84}, in which the discrete Laplacian $\laph$ is corrected where the stencil crosses $\Gamma$, taking the jumps into account. In principle we can find $u^+$ on $\Gamma$ in this manner and solve the Dirichlet problem (1.1), given $g$,  by iterating to obtain $\varphi$ so that $u^+ = g$.
W.-J. Ying et al. \cite{yh,YWjcp,YWvar,YPB} have developed a practical numerical method of this type which applies to a variety of problems in integral formulation.
An important advantage is that the integral operator is not needed explicitly, so that, for example, equations with variable coefficients can be treated.  The method we study analytically in this paper is a simplified version of the one in \cite{YWjcp}.

In its usual form the immersed interface method is uniformly $O(h^2)$ accurate
\cite{bealay,liitobook},
but in this work we use the lowest order approximation for the problem (1.3);
we correct $\laph$ for the jump in $\varphi$ itself but not for the jumps in derivatives.
Such a method was used in \cite{fed00}.  It was proved in \cite{LS} that this low order version
converges weakly to the actual solution if $\varphi$ is $C^1$.  Here we do not attempt to
approximate the exact $\varphi$ accurately; instead we find an approximate $\varphi$ to solve (1.1) accurately.  A discrete version of Green's identity based on
the technique of \cite{LS} will be important here to incorporate matching or jump conditions at
$\Gamma$ in the discrete setting.

We now describe the specific method and state the results.  The treatment is very similar for
$\R^d$ with $d = 2$ or $3$,
and we give details in $\R^2$ for simplicity with remarks about $\R^3$ as needed.
We will assume that the bounded domain $\omp$ has $C^2$ boundary $\Gamma$
and that both $\omp$ and $\omm$ are connected.   We choose a rectangular box $\bbox$ so that $\ompcl \subseteq \bbox$.
With grid size $h$ we introduce a square grid on $\bbox$.  We call the set of grid points $\boxh$ and write a grid point as $x_{i,j} = (ih,jh)$.  We choose the set $\boxh$ to include grid points on the boundary of $\bbox$; we will assign zero
boundary conditions on $\pa\bbox$.

We label each grid point in $\boxh$ as $+$ or $-$. Thus we define
$\omhp = \omp \cap \boxh$ and $\omhm = \boxh - \omp$; if a grid point is on $\Gamma$, we put it in $\omhm$.  
We partition the set $\intsh$ of intervals $[x_{i,j},x_{\ip,j}]$ and $[x_{i,j},x_{i,\jp}]$
as $\intsh = \intshp \cup \intshm \cup \intshcut$.
If interval $I$ has both endpoints in $\omhp$, then $I \in \intshp$, and similarly for $\intshm$.  If one endpoint is in $\omhp$ and the other in $\omhm$, then $I \in \intshcut$.  We call such $I$ a cut interval since it must intersect $\Gamma$. For each cut interval we find one point on $\Gamma$ in the interval;
if there are several, we simply choose one.  We call this a cut point and denote it as $x_{\iph,j}$ or $x_{i,\jph}$ to indicate it belongs to the interval.
Let $\Gamma_h^0$ be the set of all such cut points.  
It is possible that a grid point is on $\Gamma$ and is a cut point for two adjacent intervals.
Thus the number of points in $\Gamma_h^0$ could be less than the number
of intervals in $\intshcut$.
 More generally, instead of $\Gamma_h^0$, we could choose one point in each  $I \in \intshcut$ by another rule, forming a different set $\Gamma_h$.

We will solve a Poisson problem on $\boxh$ which mimics the interface problem (1.3).  Given a function $\varphi^h$
on $\Gamma_h^0$, we solve
\beq \laph u^h \eq \Phi^h \mbox{\;\;on\;} \boxh \,, \qquad u^h = 0  \mbox{\;\;on\;} \pa\boxh \eeq
where $\laph$ is the usual five-point Laplacian in $\R^2$, or seven-point in $\R^3$, and
$\Phi^h$ is determined by $\varphi^h$.  If, for example, $x_{i,j} \in \omhp$, the horizontal contribution to $h^2\laph u^h_{i,j}$
should be $u_{\ip,j}^{h,+} + u_{\im,j}^{h,+} - 2u_{i,j}^{h,+}$, where $u^{h,+}$ approximates the exact $u^+$ inside.
At $+$ points we set $u^{h,+} = u^h$, but at $-$ points we approximate $u^{h,+}$ by adding $\varphi^h$.  
Thus if $x_{\ip,j} \in \omhm$, we use
$u_{\ip,j}^{h,+} \approx u^h_{\ip,j} + \varphi^h_{\iph,j}$.  Considering
various cases similarly, we obtain the equation in (1.4) with
\beq \Phi^h_{i,j} \eq
        \sum\nolimits_{(k,\ell)} \sigma_{ijk\ell} \varphi^h_{i+k/2,j+\ell/2}\,h^{-2}  \eeq
where $(k,\ell) = (\pm 1,0)$ or $(0,\pm 1)$ and $\sigma_{ijk\ell} = 0$ if $x_{i,j}$ and $x_{i+k,j+\ell}$
have the same label; $\sigma_{ijk\ell} = -1$ if $x_{i,j}$ is $+$ and $x_{i+k,j+\ell}$ is $-$;
$\sigma_{ijk\ell} = 1$ if $x_{i,j}$ is $-$ and $x_{i+k,j+\ell}$ is $+$.  Thus
$\Phi^h_{i,j} = 0$ if the stencil of $\laph$ at $x_{i,j}$ does not cross $\Gamma$.

Having found $u^h$ on $\boxh$, we obtain values $f^h$ for the interior function $u^+$
on the set $\Gamma_h^0$ of cut points by quadratic interpolation.
Thus if $x_{i,j} \in \omhp$ and $x_{\ip,j} \in \omhm$ with $x_{\iph,j} = (ih+sh,jh)$ we define
\beq  f_{\iph,j}^h \eq u_{i,j}^h + (u_{\ip,j}^{h,+} - u_{\im,j}^{h,+})\frac{s}{2} + 
          (u_{\ip,j}^{h,+} + u_{\im,j}^{h,+}  - 2u^h_{i,j})\frac{s^2}{2}  \eeq
Here $u_{\ip,j}^{h,+} =  u_{\ip,j}^h + \varphi^h_{\iph,j}$, and $u_{\im,j}^{h,+} =  u_{\im,j}^h$ if $x_{\im,j} \in \omhp$ or $u_{\im,j}^{h,+} =  u_{\im,j}^h + \varphi^h_{\imh,j}$ if
$x_{\im,j} \in \omhm$.
With $u^{h,+}$ defined as $u^h$ on $\omhp$ and extended to the cut intervals in this way, we have  $\laph u^{h,+} = 0$
at all points in $\omhp$.  Note, however, that different values could be assigned
to $u^{h,+}$ at an endpoint in $\omhm$ by two different intervals.

To interpret $u^{h,+}$
and $f^h$, we can rewrite $\laph u^{h,+}$ near $\Gamma_h^0$, substituting for the extended values in terms of $f^h$.  If $x_{i,j} \in \omhp$ and $x_{\ip,j} \in \omhm$, with 
$x_{\iph,j} = (ih+s_rh,jh)$, and possibly $x_{\im,j} \in \omhm$,
$x_{\imh,j} = (ih-s_\ell h,jh)$, the horizontal second difference becomes 
\beq u_{\ip,j}^{h,+} + u_{\im,j}^{h,+}  - 2u_{i,j}^{h,+} \eq
\frac{2}{s_r(s_r+s_\ell)} f_{\iph,j}^h + \frac{2}{s_\ell(s_r+s_\ell)} f_{\imh,j}^h
    - \frac{2}{s_r s_\ell}u_{i,j}^h \eeq
If $x_{\im,j} \in \omhp$, $s_\ell$ is replaced by $1$ and $f_{\imh,j}$ by $u_{\im,j}^h$.
Thus $u^h$ or $u^{h,+}$, restricted to $\omhp$, is the Shortley-Weller solution of the Dirichlet problem with boundary value $f^h$ on $\Gamma_h^0$; see e.g. \cite{ciarlet,forwas,hackbook,matyam,weynans,YMSW}.  As a discretization of problems with sufficient smoothness, it is $O(h^2)$ accurate.
The equivalence of the two forms results from the fact that both use quadratic interpolation.  

Now let $A^h: \fcnscut \to \fcnscut$ be the operator $A^h\varphi^h = f^h$,
where $\fcnscut$ is the space of functions on $\Gamma^0_h$.
We prove here that $A^h$ is invertible, so that the Dirichlet problem can be solved.
We state our principal results, beginning with the last remark above.
We assume $h$ is sufficiently small; this is needed in the proof of some lemmas.  Constants below are independent of $h$.

\smallskip

{\bf Main Theorem.} {\it
(1) The grid function $u^h$ on $\omhp$, found from $\varphi^h$ according to
(1.4),(1.5), is the Shortley-Weller solution of the
discrete Dirichlet problem on $\omhp$ with boundary value $f^h$ defined in (1.6).

(2) The operator $A^h: \fcnscut \to \fcnscut$ is invertible.  Thus, given $f^h$, there is
a unique $\varphi^h$ so that $A^h\varphi^h = f^h$.

(3) There are norms $\|\cdot\|_1$ and  $\|\cdot\|_2$ on $\fcnscut$, which are discrete versions of (1.2),  such that
\beq  \|A^h\varphi^h\|_2 \leq C_1\|\varphi^h\|_1 \,,\qquad 
              \|(A^h)^{-1}f^h\|_1 \leq C_2\|f^h\|_2 \eeq

(4) To solve (1.1), assuming $g \in C^2(\Gamma)$,
let $g^h$ be the restriction to $\Gamma_h^0$.
We find $\varphi^h = (A^h)^{-1}g^h$, so that $u^h$ from
(1.4),(1.5) gives the Shortley-Weller solution of (1.1).  Then $\|g^h\|_2 \leq C_3$
and $\|(A^h)^{-1}g^h\|_1 \leq C_4$.  If $g \in C^4(\Gamma)$, $u^h$ approximates the exact $u^+$ on $\omhp$ uniformly with accuracy $O(h^2)$.
}

\smallskip

The norms $\|\cdot\|_1$, $\|\cdot\|_2$ roughly approximate the norm in the Sobolev space
$H^{1/2}(\Gamma)$, but they are not identical; they are defined in Section 3.
The boundedness of $A^h$ and its
inverse reflect the fact that $A^h$ approximates a Fredholm operator.
The $O(h^2)$ accuracy of $u^h$ is proved using the discrete maximum principle, M-matrices, or monotone matrices \cite{ciarlet,forwas,hackbook}, and more can be shown
\cite{matyam,weynans,YMSW}.
Instead of (1.6) we could use the linear interpolation
\beq  f_{\iph,j}^h \eq (1-s)u_{i,j}^h + s\,u_{\ip,j}^{h,+} \eeq
as in the method attributed to Collatz \cite{forwas,hackbook,sam}.  It also has accuracy $O(h^2)$ but is found to be less accurate \cite{sydney}.
We expect these results can be extended to the more general problem with
the equation $\Delta u = 0$ replaced by $\nabla\cdot(\beta\nabla u) = 0$, where
$\beta(x) \geq \beta_0 > 0$ is a scalar function.

The Main Theorem is proved in Sections 2--4.
In Sec. 2 we define classes of discrete functions on ${\om_h}^\pm$ with a notion of boundary value on a set $\Gamma_h^1$ different from $\Gamma_h^0$.  The new set is chosen to avoid the possibility that cut points on adjacent intervals could be very close.  A discrete Green's identity allows us to treat crude versions of the Dirichlet problem and single or double layer potentials in a Hilbert space setting.  In Lemma 2.6 we see that the discrete version of the interior integral in (1.2) is bounded by the exterior term.  The crucial importance of this fact was emphasized in \cite{cost}.
We then prove in Theorem 2.7 that a discrete version $\A$ of the operator $A$ has an inverse bounded independent of $h$ and has properties resembling those of the exact operator to a remarkable extent.  The proof is in the spirit of \cite{cost,SW}.
In Sec. 3 we relate the operator $\A$ of Theorem 2.7 to the operator $A^h$ above,
and prove
that $A^h$ is invertible, obtaining parts (2) and (3) of the Main Theorem.  In Lemma 3.1 we prove the estimate for $g^h$ in part (4).
Sec. 4 contains proofs of lemmas, including a proof that $\omhp$ is path connected for small $h$, discrete versions of the Poincar\'{e} inequality,
and the proof of Lemma 2.6.  In Sec. 5 we present numerical examples in $\R^2$ with the expected accuracy and some brief discussion.

The use of the integral equation formulation to solve a boundary value problem with the immersed interface method was suggested in \cite{mayo84}.  The augmented immersed interface method
\cite{liitobook,zliaug} is used to solve boundary value problems.  The related method of \cite{wiegbube} is also used.  The method of difference potentials \cite{diffpot,ry,hommeldiss,hommeldir} and the capacitance matrix method \cite{widcap} are similar in approach.  It was proved in \cite{hommeldir} that a difference potential version of the integral equation in $\R^2$ is solvable.  Other discrete Green's identities are given in \cite{hommeldiss,YMSW}.

\ssection{Discrete potential theory}

{\bf Discrete functions and Green's identity.}
We will use extensions of $\omhp$ and $\omhm$ past the boundary $\Gamma$.  Let $\omhpcl$ be the union of $\omhp$ with the set of endpoints of cut intervals in $\omhm$,
and similarly for $\omhmcl$.  For theoretical use we will choose a new boundary set $\Gamma_h^1$,
consisting of one point for each $I \in \intshcut$ whose distance from both endpoints is at least $\delta h$, where
$\delta$ is a chosen constant, $0 < \delta \leq \lilhalf$, independent of $h$.  We could construct $\Gamma_h^1$, for example, by choosing for $I$ the assigned point in $\Gamma_h^0$, provided it is $\delta h$ away from each endpoint, or moving it to distance $\delta h$ otherwise.  Thus $\Gamma_h^1$ is in one-to-one correspondence with $\intshcut$. This condition with $\delta > 0$
will be essential in Lemma 2.6 and Theorem 2.7, but the earlier content of this section is valid with the weaker assumption that no point in  $\Gamma_h^1$ is an endpoint.
For simplicity we will write grid functions as $f$ rather than
$f^h$ etc.

We will use extensions of functions on $\om_h^{\pm}$ to the cut intervals.  
Given a grid function $f$ on $\omhp$, for each $I \in \intshcut$ we may assign a value
of $f$ at the
endpoint of $I$ in $\omhm$.  However, this point might also be the endpoint of another interval in $\intshcut$, and we will allow a different value for each interval.  
The extended function has domain $\omhp \cup \intshcut$, but we will regard it as a function on $\omhpcl$
which is multi-valued on
$\omhpcl \cap \omhm$ and single-valued on $\omhp$.  Let $\fcnsp$ be the space of such functions.
For $f\in \fcnsp$ we can naturally define boundary values on  $\Gamma_h^1$ by linear interpolation.  In $\R^2$, if $[x_{i,j},x_{\ip,j}] \in \intshcut$ with 
$x^* = ((i+s)h,jh) \in \Gamma_h^1$,
we define the value at $x^*$ as
\beq (M_1f)_{\iph,j} \eq (1-s)f_{i,j} + sf_{\ip,j} \eeq
where the extended value is the one for this interval.
Similarly we define $(M_2f)_{i,\jph}$ for vertical cut intervals. 
For $f \in \fcnsp$, given the values on $\omhp$, the extended values determine the boundary values $M_\nu f$ and vice versa.  Thus we can think of the extensions to the endpoints as a bookkeeping device for the boundary values.  Similarly we can extend functions from
$\omhm$ to $\omhmcl$.  Let $\fcnsm$ be the space of such functions, multi-valued on
$\omhmcl \cap \omhp$, which are zero on $\pa\boxh$.  We define $M_\nu$ on $\fcnsm$.

We will use a discrete Green's identity similar to formulas in \cite{LS}.  We first introduce some notation.
For grid functions on $\boxh$ we will regard the divided differences as located on the midpoints of intervals, e.g.,
$ (D_1 u)_{\iph,j} = (u_{\ip,j} - u_{i,j})/h$, and similarly for $(D_2 u)_{i,\jph}$.  We will use the discrete
characteristic function, $\chi_{i,j} = 1$ if $x_{i,j} \in \omhp$ and $\chi_{i,j} = 0$ on $\omhm$. 
Thus $(D_1\chi)_{\iph,j} \neq 0$ iff  the interval $[x_{i,j},x_{\ip,j}]$ is in $\intshcut$.  We will need to approximate the fraction of a cut interval in $\omhpcl$.  If $I = [x_{i,j},x_{\ip,j}] \in \intshcut$, with $x^* = ((i+s)h,jh) \in \Gamma_h^1$, we define $\xi^1_{\iph,j} = s\chi_{i,j} + (1-s)\chi_{\ip,j}$ and similarly we define $\xi^2_{i,\jph}$.  If $I \in \intshp$, $\xi^\nu = 1$, while if $I \in \intshm$, $\xi^\nu = 0$.
We note for later use (cf. pp. 1738-39 in \cite{LS}) that 
\beq D_1(f\chi)_{\iph,j} \eq (D_1f)_{\iph,j} \xi^1_{\iph,j} \p
     (M_1f)_{\iph,j} (D_1\chi)_{\iph,j} \eeq
and similarly for $D_2$.  This description given in $\R^2$ extends directly to $\R^3$.

The discrete Green's identity in $\R^d$, valid for functions $u, v$ on $\boxh$ and also for multi-valued functions in $\fcnsp$, is, with $1 \leq \nu \leq d$,
\beq   \sum_{i,j} (\laph u) v\chi \,h^d \p 
  \sum_{i,j, \,\nu} (D_\nu u)(D_\nu v)\xi^\nu \,h^d \eq
\m \sum_{i,j,\,\nu} (D_\nu u) (M_\nu v) (D_\nu \chi) \,h^d  \eeq
The first sum is over grid points in $\omhp$ and the second is on intervals in
$\intshp\cup\intshcut$.  Since $D_\nu\chi = \pm 1/h$ on the cut intervals and zero otherwise, the
right side is a sum over $\intshcut$ and is the analog of a boundary integral.  (The sign is opposite from usual because $\chi$ decreases going outward.)
The complementary identity for $u,v \in \fcnsm$ is 
\beq
  \sum_{i,j} (\laph u) v(1 - \chi) \,h^d
    \p \sum_{i,j, \,\nu} (D_\nu u)(D_\nu v)(1 - \xi^\nu) \,h^d
        \eq \sum_{i,j, \,\nu} (D_\nu u) (M_\nu v) (D_\nu \chi) \,h^d 
\eeq

We verify (2.3) in $\R^2$, assuming at first the functions are single-valued.
We start with the identity, for fixed $j$,
\beq \sum\nolimits_i \,(D_1^2 u)_{i,j}v_{i,j} \eq 
      \m \sum\nolimits_i \, (D_1 u)_{\iph,j}(D_1 v)_{\iph,j} \eeq
assuming $ v = 0$ on $\pa\boxh$, where $D_1^2$ is the usual centered second difference on $\boxh$.  We replace $v$  with $\chi v$ and use the identity (2.2).  We then sum over $j$, obtaining an expression like (2.3) with $D_1^2 u$ on the left and
$\nu = 1$. We repeat for $\nu = 2$, sum over $\nu$, and set $\laph = \Sigma_\nu D_\nu^2$ to obtain (2.3). For a multi-valued function in dimension one, we note that (2.3) holds for an interval of the form
$[x_1,\dots,x_k]$ if $x_1 \in \omhm$, $x_k \in \omhm$, and $x_i \in \omhp$ for
$1 < i < k$, and the sum over the full interval is the sum over a union of such intervals.  For the general case we can repeat the other sums successively as before.

Let $\fcnsgam$ be the space of functions on $\Gamma_h^1$.  Given $f \in \fcnsp$ we have the boundary value $Mf \in \fcnsgam$, as defined in (2.1).  Conversely, given $\psi \in \fcnsgam$, there
are $f \in \fcnsp$ with $Mf = \psi$.  For example we can set $f = 0$ on $\omhp$ and then choose $f$ on $\omhpcl \cap \omhm$, allowing multi-values, so that $Mf = \psi$.  Similarly there exist $f \in \fcnsm$ with $Mf = \psi$.

\medskip

{\bf Hilbert space structure and the Dirichlet problem.}
The identities (2.3), (2.4) suggest using discrete Dirichlet norms as inner products.  To proceed we need the following lemmas.  Recall that we assume $\omp$ and $\omm$ are connected and $\Gamma$ is $C^2$.  
The lemmas require $h$ to be sufficiently small depending on the geometry
of the domains $\omp$ and$\omm$, as will be seen in the proofs.

\begin{lemma}[\bf Connectedness]   
For $h$ sufficiently small, the grid points in $\omhp$ and $\omhm$ are path connected:  Given two grid points in $\omhp$, there is a path between them consisting of grid intervals with all endpoints in $\omhp$, and similarly for $\omhm$.
\end{lemma}

\begin{corollary}  
For $h$ small, if $u$ is a function on $\omhp$ with $D_\nu u \equiv 0$ on $\intshp$, $1 \leq \nu \leq d$,
then $u$ is constant, and similarly for $\omhm$.
\end{corollary}

\begin{lemma}[\bf Poincar\'{e} inequality]
There is a constant $C$ independent of $h$ so that for $h$ sufficiently small
the following hold.

(1) For any $u \in \fcnsp$
with $Mu = 0$, we have
\beq \sum\nolimits_{\,\omhp} \, u^2 \,h^d \lee C \sum (D_\nu u)^2 \xi^\nu\,h^d \eeq

(2) For any function $u$  on $\omhp$ with mean value zero, i.e., $\sum u = 0$ on $\omhp$, we have
\beq \sum\nolimits_{\,\omhp} \, u^2 \,h^d \lee C \sum\nolimits_{\,\intshp}\, (D_\nu u)^2 \,h^d \eeq

(3) For any function $u$ on $\omhm$ with $u = 0$ on $\pa\boxh$ the inequality (2.7) holds with sums on $\omhm$ and $\intshm$.
\end{lemma}

The lemmas are proved in Sec. 4, except that we briefly prove (1) in Lemma 2.3.  Given
$u \in \fcnsp$, define $\utw$ on $\boxh$ as $\utw = u$ on $\omhp$ and $\utw = 0$ otherwise.
On a cut interval, using $M_\nu u = 0$, we find $D_\nu \utw = \xi^\nu D_\nu u$, so that
$\sum_{\boxh}( D_\nu \utw)^2\,h^d \leq C \sum_{\intshp} (D_\nu u)^2 \xi^\nu \,h^d$.  Since $\utw = 0$ near
$\pa\boxh$, it is a standard fact that $\sum_{\boxh}\utw^2\,h^d \leq \sum_{\boxh}( D_\nu \utw)^2\,h^d$, and the conclusion follows.

Assuming $h$ is small enough for the lemmas to hold,
we now define inner products and seminorms, for $u^+, v^+ \in \fcnsp$ and $u^-, v^- \in \fcnsm$,
\begin{eqnarray}
 \langle u^+,v^+ \rangle^+ \eq \sum (D_\nu u^+)(D_\nu v^+)\xi^\nu  \,h^d \qquad \\
 \langle u^-,v^- \rangle^- \eq \sum (D_\nu u^-)(D_\nu v^-)(1-\xi^\nu) \,h^d  
\end{eqnarray}
with $\|u^+\|^+ = [\langle u^+,u^+ \rangle^+]^{1/2}$ and similarly for $\|u^-\|^-$.
Then $\|\cdot\|^-$ is a norm, by Lemma 2.3(3), while $\|\cdot\|^+$ is a seminorm.
Let $\fcnspz$ be the subspace of $\fcnsp$ with $Mf = 0$.  Cor. 2.2 implies that 
$\|\cdot\|^+$ is a norm on this subspace:  If $D_\nu f \equiv 0$, then 
$f$ is constant, but $Mf = 0$ so that the constant is $0$.  

We can now treat the Dirichlet problem  for discrete harmonic functions
on $\omhp$ or $\omhm$.  Given
$\psi \in \fcnsgam$, the two problems are to find $v^+ \in \fcnsp$ and $v^- \in \fcnsm$ so that
\begin{eqnarray}
\laph v^+ \eq 0 \mbox{\;\;on\;} \omhp \,,\qquad  Mv^+ = \psi \mbox{\;\;on\;} \Gamma_h^1 
\qquad\qquad\\
\laph v^- \eq 0 \mbox{\;\;on\;} \omhm \,,\qquad  Mv^- = \psi \mbox{\;\;on\;} \Gamma_h^1 \,,    \qquad v^- = 0  \mbox{\;\;on\;} \pa\boxh 
\end{eqnarray}

\begin{lemma}
Given $\psi \in \fcnsgam$, each of the problems (2.10),(2.11) has a unique solution.
\end{lemma}

Proof:  We can adapt the standard proof to this context.  For the problem in $\omhp$, we have seen there is some $w \in \fcnsp$ with $Mw = \psi$.  Setting $v^+ = v_0 + w$, our problem is to find
$v_0 \in \fcnspz$ such that $\laph v_0 = - \laph w$ on $\omhp$.  We seek $v_0$ in the weak form
$$  \sum (D_\nu v_0)(D_\nu z)\xi^\nu \,h^d \eq \sum (\laph w) z\chi \,h^d \,,
         \qquad z \in \fcnspz        $$
The existence of such $v_0$ follows from the fact that the left side is the inner product in the Hilbert space $\fcnspz$.  Then from (2.3) we have
$$ \sum \laph (v_0 + w) z\chi \,h^d \eq 0 \,, \qquad \mbox{for all\;} z \in \fcnspz $$
We need to check $f \equiv \laph (v_0 + w)$ is zero on $\omhp$.  If $x \in \omhp$ is not the endpoint of a cut interval, we define $z = 1$ at $x$ and $z = 0$ elsewhere; we conclude
$f(x) = 0$.  For other $x \in \omhp$, set $z(x) = 1$ and, for any cut interval with $x$ as an endpoint, define $z$ at the endpoint in $\omhm$ so that $Mz = 0$ for that interval; set $z = 0$ otherwise. Again we find $f(x) = 0$, and the claim is verified.  Uniqueness is evident from the construction.
The exterior problem is similar.

\medskip

The discrete harmonic functions have a minimizing property like that of the exact harmonic functions.

\begin{lemma}
Suppose $v, w \in \fcnsp$, $\laph v = 0$ on $\omhp$,
and $Mv = Mw$ in $\fcnsgam$.  Then
$\langle v,v \rangle^+  \leq  \langle w,w \rangle^+$.  A similar statement holds in $\fcnsm$.
\end{lemma}

Proof:  Setting $w_0 = w - v$, we have $Mw_0 = 0$, so that $\langle v,w_0 \rangle^+ = 0$ from Green's identity (2.3),
and thus $\langle w,w \rangle^+ = \langle v,v \rangle^+ + \langle w_0,w_0 \rangle^+$.

\medskip

Given $\psi \in \fcnsgam$ the harmonic pair $v^+, v^-$ is determined so that $Mv^+ = Mv^- = \psi$. 
We can think of $v^\pm$ as a crude version of a single layer potential.  We will write $v^\pm = S\psi$,
meaning $v^\pm$ is the single layer determined by $\psi$ as in (2.10),(2.11), so that $Mv^\pm = \psi$. 
We will also write $v^+ = S^+\psi$, $v^- = S^-\psi$.  
Note that for $v^\pm \in \fcnspm$ with $\laph v^\pm = 0$ we have 
$S^\pm M v^\pm  = v^\pm$ by uniqueness.

Motivated by \cite{cost,SW}
we define an inner product on $\fcnsgam$, making it a Hilbert space.  If $v^\pm = S\psi$ and $z^\pm = S\zeta$,
the inner product $(\psi,\zeta)$ is that of the single layers,
\begin{multline}
 (\psi,\zeta)  \eq \langle v^+,z^+ \rangle^+ + \langle v^-,z^- \rangle^- \\
   \eq \sum (D_\nu v^+)(D_\nu z^+)\xi^\nu \,h^d + \sum (D_\nu v^-)(D_\nu z^-)(1-\xi^\nu) \,h^d 
\end{multline}
The corresponding norm will be important to us,
\beq \|\psi\|^2 \eq \sum (D_\nu v^+)^2\xi^\nu \,h^d + \sum (D_\nu v^-)^2 (1-\xi^\nu)\,h^d \,,
       \qquad   v^\pm = S\psi
\eeq
From (2.3),(2.4) the inner product reduces to a discrete boundary integral,
\beq (\psi,\zeta)  \eq  -\sum [D_\nu v] (M_\nu z) (D_\nu \chi) \,h^d
                \eq  -\sum [D_\nu v] (\zeta) (D_\nu \chi) \,h^d  \eeq
where $[D_\nu v] = D_\nu v^+ - D_\nu v^-$ and similarly with $\psi,\zeta$ reversed.

\medskip

{\bf The discrete double layer potential and invertibility.}
Given $\varphi \in \fcnsgam$ we can solve the interface problem (1.4),(1.5).
(Note with $\Gamma_h^1$, we do not have the identification of cut points common to two intervals imposed in Sec. 1.)
From the solution $u^h$ we can define
$u^+ = u^h$ on $\omhp$ and extend to $\omhpcl$ by adding $\varphi$ at the extended endpoints:
For a cut interval $[x_{i,j},x_{\ip,j}]$ in $\R^2$
with $x_{ij} \in \omhp$ and $x_{\ip,j} \in \omhm$ we define $u^+_{\ip,j} = u^h_{\ip,j} + \varphi_{\iph,j}$ and similarly
for other intervals, so that $u^+$ is defined in $\fcnsp$.  We define $u^- \in \fcnsm$ in the same manner,
subtracting $\varphi$ at endpoints of cut intervals in
$\omhp$.  From the construction of $u^h$ we have $\laph u^+ = 0$ on $\omhp$ and $\laph u^- = 0$ on $\omhm$.
For the interval described, we have $u^-_{i,j} = u^h_{i,j} - \varphi_{\iph,j}$, or
$u^+_{i,j} - u^-_{i,j} = \varphi_{\iph,j}$, and 
$u^+_{\ip,j} = u^-_{\ip,j} + \varphi_{\iph,j}$; the same holds generally. It follows that $Mu^+ - Mu^- = \varphi$
and $Du^+ = Du^-$ on the cut intervals. In summary, we have $u^+ \in \fcnsp$ and $u^- \in \fcnsm$ so that
\begin{eqnarray}
\laph u^+ \eq 0 \mbox{\;\;on\;} \omhp \,,\qquad  \laph u^- \eq 0 \mbox{\;\;on\;} \omhm \,,
\qquad\qquad\nonumber \\
 \, [Mu] \eq \varphi \mbox{\;\;on\;} \Gamma_h^1 \,,\qquad [Du] \eq 0 \mbox{\;\;on\;} \intscut \,,
       \qquad u^- = 0  \mbox{\;\;on\;} \pa\boxh 
\end{eqnarray}
It is important that $u^+$, $u^-$ do not depend on the choice of $\Gamma_h^1$.
Conversely, if $u^+, u^-$ satisfy (2.15) then we can define $u^h = u^+$ on $\omhp$ and $u^h = u^-$ on $\omhm$
to recover the original $u^h$.  Thus the two formulations are equivalent. We see that $u^+$, $u^-$ form an
analog of the double layer potential. 
From Green's identities (2.3),(2.4) we have
\beq \langle u^+,u^+ \rangle^+ + \langle u^-,u^- \rangle^- \eq 
             \m \sum (D_\nu u^\pm) \varphi (D_\nu \chi) \,h^d \eeq
where $D_\nu u^\pm = D_\nu u^+ = D_\nu u^-$.  It follows that $u^+,u^-$ are uniquely determined by $\varphi$;
if $\varphi = 0$, then $\langle u^-,u^- \rangle^- = 0$ implies $u^- = 0$, but then $Mu^+ = 0$, so that $u^+ = 0$ as well.

We now define operators $\A, \B : \fcnsgam \to \fcnsgam$, as
$\A(\varphi) = Mu^+$ and $\B(\varphi) = Mu^-$, with $u^\pm$ as in (2.15).  They are the analogues of the classical Fredholm operators.  Note ${\cal A}-{\cal B} = I$ and $(\A + \B)/2$ corresponds to the classical integral operator.
We prove that ${\cal A}$ is invertible and has properties like those in the classical case.
The following lemma will be fundamental.  It is proved in Sec. 4.  The importance of this fact in classical potential theory
was emphasized in \cite{cost}; see Lemma 1.

\begin{lemma}[\bf Extension Lemma]
There is a constant $C_1$, independent of $h$,
so that for any $\varphi \in \fcnsgam$,
and the corresponding single layer potential $v^+ \in \fcnsp$ and $v^- \in \fcnsm$, such that $\laph v^+ = 0$ on $\omhp$,
$\laph v^- = 0$ on $\omhp$, and $Mv^+ = Mv^- = \varphi$, we have
\beq \langle v^+,v^+ \rangle^+ \leq C_1 \langle v^-,v^- \rangle^-    \eeq
and if $v^+$ has mean value zero on $\omhp$,
\beq \langle v^-,v^- \rangle^- \leq C_1 \langle v^+,v^+ \rangle^+    \eeq
\end{lemma}

\begin{theorem}[\bf Bounded invertibility of $\A$]
The operators $\A$ and $\B$ on $\fcnsgam$ with inner product (2.12) are symmetric.  In the operator norm (2.13),
$\|\B\| \leq r$ where $r < 1$ is independent of $h$, so that $\A = I + \B$ is invertible with $\|\A^{-1}\|\leq (1-r)^{-1}$.
The spectrum of $\A$ is in the interval $1-r \leq \lambda \leq 1$.
\end{theorem}

Proof.  The proof is inspired by \cite{cost,SW}.  Given $\varphi \in \fcnsgam$, we will construct the double layer $u^\pm$ as in (5.3) in \cite{SW}.
Let $v^\pm = S\varphi$, as in (2.10),(2.11), so that $Mv^\pm = \varphi$.  We wish to find $w^\pm$, a
single layer of the form $w^\pm = S\omega$ for some $\omega \in \fcnsgam$, such that $[Dw] = -Dv^+$ on $\intshcut$.
In view of (2.14), we seek $\omega$ so that for all $\zeta \in \fcnsgam$,
\beq (\omega,\zeta) \eq \sum (D_\nu v^+) \zeta (D_\nu \chi) \,h^d \eeq
Such $\omega$ exists since the right side is a linear functional on $\zeta \in \fcnsgam$.  Then from (2.14) we have
$$ \sum \left([D_\nu w] +(D_\nu v^+)\right) \zeta (D_\nu \chi)\,h^d = 0$$
Since $\zeta$ is arbitrary, we conclude $[Dw] = - Dv^+$, verifying the choice of $\omega$ and $w^\pm$. 
We now set $u^- = w^-$ and $u^+ = v^+ + w^+$, so that $Mu^+ - Mu^- = Mv^+ + (Mw^+ - Mw^-) = \varphi$,
and $Du^+ - Du^- = Dv^+ + [Dw] = 0$, and by uniqueness $u^\pm$ is the double layer $\A\varphi$ as in (2.15).

Note that $\B\varphi = \omega$, since $u^- = w^-$ and thus $Mu^- = Mw^- \eq \omega$.  We will show that $\B$ is a contraction.
With $z^\pm = S\zeta$, so that $Mz^\pm = \zeta$, we use (2.3) to write (2.19) as
\beq (\omega,\zeta) \eq - \sum (D_\nu v^+)(D_\nu z^+) \chi\,h^d 
             \eq - \langle v^+,z^+ \rangle^+  \eeq
We set $\zeta = \omega$, so that $\|\omega\|^2 = -\langle v^+,w^+ \rangle^+$,
and recalling the definition (2.13) of $\|\omega\|$ from $w^\pm$, we find that
$ \|\omega\|^2 \leq \langle v^+,v^+ \rangle^+$.  
We note that $\|\varphi\|^2 = \langle v^+,v^+ \rangle^+ + \langle v^-,v^- \rangle^-$; call the two terms $m^+$ and $m^-$.  From Lemma 2.6 we have $m^+ \leq C_1 m^-$.  It follows that $m^+/(m^- + m^+) \leq C_1/(1 + C_1)$
and $m^+ \leq r^2 \|\varphi\|^2$ with $r^2 = C_1/(1 + C_1) < 1$.  Thus $\|\omega\| \leq r\|\varphi\|$.
Since $\omega = \B\varphi$, we have proved that $\|\B\| \leq r < 1$.  Moreover, we can interpret the identity (2.20)
as $(\B\varphi,\zeta) =  - \langle (S\varphi)^+,(S\zeta)^+ \rangle^+$, from which it is evident that $\B$ is symmetric
and $\B \leq 0$.  

Since $\A = I + \B$, the properties of $\B$ imply that $\A$ is bounded and symmetric,
$\A$ is invertible, and $\|\A^{-1}\| \leq (1-r)^{-1}$.  The spectrum of $\B$ is limited to
$-r \leq \lambda \leq 0$, and thus the spectrum of $\A$ is in the interval
$1-r \leq \lambda \leq 1$.

\ssection{Proof of the Main Theorem}

To prove the Main Theorem we will relate the operator $A^h: \fcnscut \to \fcnscut$
of Sec. 1 to
$\A: \fcnsgam \to \fcnsgam$ of Sec. 2 since the latter is invertible.  We use the notation of Sec. 2.  We will denote the norm on $\fcnsgam$,
defined by (2.13),(2.10),(2.11) as $\|\cdot\|_{SL}$
to distinguish from other norms.  We will write $A$ rather than $A^h$ etc.

In applying the theory of Sec. 2 to the method of Sec. 1 we need to regard boundary functions as defined on $\Gamma_h^0$ rather than $\Gamma_h^1$ or $\intshcut$.  We noted that a point in $\Gamma_h^0$ could be a multiple cut point, i.e., a cut point for more than one interval in $\intshcut$.  Such a point must be a grid point in $\omhm$.  For this reason we define the subspace $\fcnseq$ of $\fcnsgam$, consisting of functions $\psi$ such that, if $x_p, x_q \in \Gamma_h^1$ belong to intervals with the same cut point in $\Gamma_h^0$, then $\psi(x_p) = \psi(x_q)$.  We can identify
$\fcnscut$ with $\fcnseq$ in a natural way; they have the same dimension.

We will determine the image subspace $\A(\fcnseq)$.  Suppose $\psi \in \fcnsgam$, with double layer potential $u^+ \in \fcnsp$, $u^- \in \fcnsm$ as in (2.15).
If $x \in \Gamma_h^0 \cap \omhm$ is a multiple cut point, and thus a grid point, then $u^-(x)$ is single-valued, and the multiple values of $u^+(x)$ have the form $u^-(x) + \psi$.  If $\psi \in \fcnseq$, these values are equal, and $u^+(x)$ is single-valued.  Let $\fcnspeq$ be the subspace of $\fcnsp$ with this property.  Conversely, if $u^+$ is single-valued at the multiple cut points, it follows that $\psi \in \fcnseq$.  Since $\A\psi = Mu^+$ and $Mu^+$ determines $u^+ = S^+(Mu^+)$, we have identified $\A(\fcnseq)$ as the elements of $\fcnsgam$ whose harmonic extensions to $\fcnsp$ are in $\fcnspeq$. 

Given $\varphi \in \fcnscut$, we can express $A\varphi \in \fcnscut$ using $\A$.  Form $\phitw \in \fcnseq$ from $\varphi$ by repeating values as needed for multiple cut points.  Let $u^\pm$ be the double layer potential determined by $\phitw$, so that $u^+ \in \fcnspeq$.  Let $Q: \fcnspeq \to \fcnscut$ be the quadratic interpolation defined by (1.6).
Then $A\varphi = Qu^+ = QS^+(\A\phitw)$.

Next we define the norms to be used.  For $\varphi \in \fcnscut$ we define
$\|\varphi\|_1 = \|\phitw\|_{SL}$, with $\phitw$ as above.  
For the second norm, suppose
$f \in \fcnscut$; let $w$ be the Shortley-Weller solution of the discrete Dirichlet problem
on $\omhp$ with boundary value $f$.  We extend $w$ to an element of $\fcnspeq$ by quadratic
extrapolation:  If $x_{i,j} \in \omhp$, $x_{\ip,j} \in \omhm$ and $x_{\im,j} \in \omhm$,
we define $w$ at $x_{\ip,j}, x_{\im,j}$ from the quadratic function determined by 
$f_{\iph,j}, f_{\imh,j}$, $w_{i,j}$.  In the more usual case $x_{i,j} \in \omhp$, $x_{\ip,j} \in \omhm$ and $x_{\im,j} \in \omhp$, $w_{\ip,j}$ is defined from $f_{\iph,j}$, $w_{i,j}, w_{\im,j}$.
We note for later that $Qw = f$ because of the uniqueness of the quadratic fit.  The extended
function has $\laph w = 0$ on $\omhp$ since the two forms of the second difference are equivalent.
Now with $Mw$ defined as in (2.1), we set $\|f\|_2 = \|Mw\|_{SL}$.

We check that $\|A\varphi\|_2 \leq C \|\varphi\|_1$, i.e., $A$ is bounded in the appropriate sense:  If $A\varphi = Qu^+$ as above, and $f = Qu^+$, then the Shortley-Weller solution for
$f$ is $w = u^+$, so that
$\|Qu^+\|_2 = \|Mu^+\|_{SL} = \|\A\phitw\|_{SL} \leq 2\|\phitw\|_{SL} = 2\|\varphi\|_1$, where we use the boundedness of $\A$ in the inequality.

Next we check that $A$ is invertible and
$\|A^{-1}f\|_1 \leq C\|f\|_2$.  Given $f \in \fcnscut$, we seek $\varphi \in \fcnscut$ so that $A\varphi = f$.  Let $w \in \fcnspeq$ be the extended Shortley-Weller solution with boundary value $f$, so that $Qw = f$.  Let $\phitw = \A^{-1}(Mw)$.
Then $\A\phitw = Mu^+$ where $u^\pm$ is the double layer determined by $\phitw$.
Thus $Mu^+ = \A\phitw = Mw$, and therefore $u^+ = w \in \fcnspeq$ by uniqueness of the 
Dirichlet problem in $\fcnsp$, Lemma 2.4.
Since $u^+ \in \fcnspeq$, we have $\phitw \in \fcnseq$
as noted above, and thus there is $\varphi \in \fcnscut$ corresponding to $\phitw$. 
In summary, $A\varphi = Qu^+ = Qw = f$, as required.  As for the estimate, we know
that $\|\phitw\|_{SL} \leq C \|\A\phitw\|_{SL}$ since $\|\A^{-1}\|$ is bounded.  
With $\|f\|_2 = \|Mw\|_{SL} = \|\A\phitw\|_{SL}$ and 
$\|\varphi\|_1 = \|\phitw\|_{SL}$, we conclude that $\|\varphi\|_1 \leq C\|f\|_2$.

We have now verified (2) and (3) in the Main Theorem.  The next lemma proves (4).

\medskip

\begin{lemma}
Suppose $\fbar$ is a $C^2$ function on $\Gamma$. With $h$ chosen, let $f$ be the
restriction of $\fbar$ to $\Gamma_h^0$.  Let $w \in \fcnsp$ be the Shortley-Weller solution
on $\omhp$, extended to $\omhpcl$ by quadratic extrapolation. Then $\|f\|_2 = \|Mw\|_{SL}$
is bounded independent of $h$.
\end{lemma}

Proof.  It is well known that $w$ is uniformly bounded on $\omhp$, independent of $h$; this
is shown using the maximum principle or the theory of monotone and M-matrices
\cite{ciarlet,forwas,hackbook}.
We can assume $\fbar$ has a $C^2$ extension to a neighborhood of $\ompcl$.  We write
$w = v + \fbar$, so that $v$ on $\omhp$
is the Shortley-Weller solution of $L_h v = -L_h \fbar$ with
zero boundary value, where $L_h$ is $\laph$, modified at the boundary.  Now let $\vtw = v$ on $\omhp$ and $\vtw = 0$ on $\omhm$.
It is proved in \cite{hackbook}, pp. 294-5, that
\beq \sum\nolimits_{\boxh} (D_\nu \vtw)^2\, h^d \leq 
\left| \sum\nolimits_{\boxh} (L_h \fbar)\vtw\,h^d \right| \eeq
which is bounded.  
In particular 
\beq \sum\nolimits_{\intshp} (D_\nu w)^2\, h^d \lee C \eeq

We need to extend this estimate to the cut intervals.  We use the fact that
$\laph w = 0$ on $\omhp$ to relate these differences to those in (3.2).
We consider cases in $\R^2$:  (1) Suppose $x_{i,j} \in \omhp$,
$x_{\ip,j} \in \omhm$, and the rest of the stencil of $\laph w_{i,j}$ is in $\omhp$.
Thus there is one cut point, with the form $x_{\iph,j} = (ih+sh, jh)$.
Let $d_1 = D_1w_{\iph,j} = (w_{\ip,j} - w_{i,j})/h$.  Then $d_1 = - d_2 - d_3 - d_4$, where
each $d_k$ on the right is a difference appearing in the sum (3.2).  

(2) Suppose $\laph w_{i,j}$ 
has two cut intervals, in different directions,
such as $x_{\iph,j} = (ih+s_1h,jh)$ and
$x_{i,\jph} = (ih,jh+s_2h)$, where $s_1, s_2 > 0$ could be arbitrarily small.
The quadratic interpolation gives 
$D_1w_{\iph,j} = a_1(f_{\iph,j} - w_{i,j})/h + b_1 D_1w_{\imh,j}$ where $a_1 = 2/(s_1(1+s_1))$
and $b_1 = (s_1-1)/(s_1+1)$.  Since $|f_{\iph,j} - f_{i,j}| \leq C s_1h$, we can replace the
first difference by $a_1(f_{i,j} - w_{i,j})/h$ plus a bounded term.  Treating $D_2w_{i,\jph}$
similarly and writing $0 = \laph w = \sum_k d_k$ as before, we conclude that $(a_1 + a_2)(f_{i,j} - w_{i,j})/h$ is a sum of terms which are bounded or differences in (3.2).  Since $a_1,a_2 > 0$, the same is true for $D_1w_{\iph,j}$ and
$D_2w_{i,\jph}$.  

(3) Suppose $x_{i,j} \in \omhp$, $x_{\ip,j} \in \omhm$  $x_{\im,j} \in \omhm$.
Then assuming $h$ is small, $\Gamma$ is close to horizontal, and either $x_{i,\jp} \in \omhm$ and $x_{i,\jm} \in \omhp$ or the reverse;
we assume the former, with $x_{i,\jph} = (ih,jh+s_2h)$.  Since we are estimating differences, we can assume $w_{i,j} = 0$.
With $x_{\iph,j} = (ih+s_r h,jh)$, $x_{\imh,j} = (ih-s_\ell h,jh)$, we find the extrapolated values $w_{\ip,j} = a_r f_0 + f'_0 h + b_{\ip,j}h^2$,
$w_{\im,j} = a_\ell f_0 - f'_0 h + b_{\im,j}h^2$, where $f_0 = f_{i,j}$, $f'_0$ is the horizontal derivative of $\fbar$ at $x_{i,j}$, $a_r = (1+s_\ell - s_r)/(s_rs_\ell)$,
$a_\ell = (1+s_r - s_\ell)/(s_r s_\ell)$, and $|b_{i \pm 1,j}|$ is bounded in terms of
the second derivative of $\fbar$.  Treating $w_{i,\jp}$ as before, we find that $\laph w_{ij}$ has
terms $(a_r + a_\ell + a_2)f_0/h$ and terms which are bounded or differences in (3.2),
with $a_r,a_\ell,a_2 > 0$.  Proceeding as in case (2), we conclude that the three
differences on cut intervals are bounded. 
 
A similar argument can be
used in $\R^3$; for small $h$, the stencil of $\laph$ always has one interval which is not cut.
In summary we have shown that, for each $I \in \intshcut$, the difference is estimated as
$|D_\nu w|^2 \leq C_1\sum |D_\nu w|^2 + C_2$, with the sum over intervals in $\intshp$ connected to $I$.  Summing over $I \in \intshcut$ and using (3.2) we can conclude that
\beq   \langle w,w \rangle^+  \eq \sum (D_\nu w)^2 \xi^\nu \, h^d  \leq C  \eeq

It remains to estimate the extension of $w$ to $\omhmcl$.  We set $w = \wtw + w_0$ in $\fcnsp$, where $w_0$ is the mean value of $w$ on $\omhp$.  We can apply Lemma 2.6 to $\wtw$.  Let $\wtw^- \in \fcnsm$ be the solution of $\laph \wtw^- = 0$ in $\omhm$ with $M\wtw^- = M\wtw$.  Then 
$\langle \wtw^-,\wtw^-\rangle^-$ is bounded by (3.3), according to (2.18).  Finally, $w_0$ is bounded, and its extension to
$\fcnsm$ is $w_0 v^-$, where $v^- \in \fcnsm$, $\laph v^- = 0$, and $Mv^- = 1$.  We can check that
$\langle v^-, v^- \rangle^-$ is bounded by defining a smooth function $\zeta$
which is $\equiv 1$ in a neighborhood of
$\Gamma$ and $\equiv 0$ in a neighborhood of $\pa\bbox$, applying Green's identity (2.4) to
$v^- - \zeta$, and using Lemma 2.3(3) to estimate $\|v^- - \zeta\|^-$.
We now have $w^- = \wtw^- + w_0 v^-$ bounded with
$Mw^- = Mw$, and the proof is complete.

\ssection{Proofs of the Lemmas}

We begin with an observation that will be used in the proofs of Lemmas 2.1 and 2.3,
first in $\R^2$.
We have assumed $\Gamma$ is $C^2$.  We note that
each point $z \in \Gamma$ has an open neighborhood
${\cal G}_z$ where $x = (x^1,x^2) \in \Gamma$ if $x^2 = \gamma(x^1)$,
$x \in \omp$ if $x^2 < \gamma(x^1)$ and $x \in \omm$ if $x^2 > \gamma(x^1)$,
with some function $\gamma$,
or a similar form with $x^1, x^2$ or the inequalities reversed.
We can assume the ${\cal G}_z$ are small enough so that the mapping along normal lines is
invertible, and we assume ${\cal G}_z$ is a union of normal lines.  Then we define
${\cal N}_z = \{(x^1,\gamma(x^1)-y): |x^1-z^1|<\alpha_z, |y| <\beta_z\}$ with
$\alpha_z, \beta_z > 0$ chosen small enough so that ${\cal N}_z \subseteq {\cal G}_z$.
We can cover $\Gamma$ by the union $S$
of finitely many of the ${\cal N}_z$.
We will show that there is a constant
$C_0$, depending only on the maximum $m$ of the various $|\gamma'|$, so that
\beq  y \leq C_0\, \mbox{dist}(x,\Gamma) \quad \;
         \mbox{for\;} \;\;  x = (x^1,\gamma(x^1) - y) \in {\cal N}_z\cap\omp  \eeq
where $\mbox{dist}$ is the distance.
Since $S$ is open we can choose $\rho > 0$ such that $x \in S$
provided $\mbox{dist}(x,\Gamma) \leq \rho$.  We also require $\rho$ small enough so that $2C_0\rho < \beta_z$ for the finitely many chosen $z$.  The same construction can be done in $\R^3$.

To verify (4.1), suppose
$(x^1,\gamma(x^1)-y) \in {\cal N}_z$ with $y > 0$ and
$(x^0,\gamma(x^0)) \in \Gamma\cap{\cal G}_z$.  Let $r$ be the distance.  Then
$r^2 = (x^1 - x^0)^2 + (y-s)^2$, where $s = \gamma(x^1) - \gamma(x^0)$.  With $\sig > 0$
arbitrary,
$2y|s| \leq (1+\sig)^{-1}y^2 + (1+\sig)s^2$, so that
$r^2 \geq (x^1 - x^0)^2 + (1-(1+\sig)^{-1})y^2 - \sig s^2$.  Now 
$|s| \leq m|x^1-x^0|$, and choosing  $\sig = 1/m^2$
we have $r^2 \geq \sig(1+\sig)^{-1}y^2$, which is equivalent to $y \leq C_0 r$ with
$C_0 = \sqrt{m^2 + 1}$.

\bigskip

{\bf Proof of Lemma 2.1.}
We give the proof for $\omp \subseteq \R^d$.  Since $\omp$ is open and connected, it is path connected.
Let $\om_\rho = \{x \in \omp: \mbox{dist}(x,\Gamma) > \rho\}$.  We can map $\om_\rho$ bijectively
to $\omp$, and thus it is also path connected.
We choose $\eps$ with $0 < \eps < \rho/(3\sqrt{d})$.
Since ${\overline {\om_\rho}}$ is compact, we can choose a finite set
$Y$ of points $y_j \in {\overline {\om_\rho}}$ so that the union of the balls of radius $\eps$
about the $y_j$'s cover ${\overline {\om_\rho}}$.  The choice of $\eps$ ensures that  
the square box $B_j$ centered at $y_j$ with side
$6\eps$ is  $\subseteq \omp$ for each $j$.

We first show that any two grid points in $\om_\rho$ are connected by a grid path in $\omhp$.
Suppose we are given $x_0$ and $z$ in $\om_\rho \cap \omhp$.
Then there is a path $x(t)$ in $\om_\rho$,
$0 \leq t \leq 1$ such that $x(0) = x_0$ and $x(1) = z$.
By uniform continuity there exists $\tau > 0$ so that
$|t - t'| \leq \tau$ implies $|x(t) - x(t')| < \eps$.
We know there is some $y_0 \in Y$ so that $|y_0 - x_0| < \eps$.  Consider
$t = \tau$; 
$|x(\tau) - x(0)| \leq \eps$, so that $|x(\tau) - y_0| < 2\eps$.  We can pick a grid point $\xtw_1 \in \omhp$ so that $|\xtw_1 - x(\tau)| < \eps$, assuming $h$ small enough
so that $\sqrt{d}h/2 < \eps$.  Then
$|\xtw_1 - y_0| < \eps + |x(\tau) - y_0| < 3\eps$.
Thus $x_0 = x(0)$ and $\xtw_1$ are both in the ball of radius $3\eps$ about $y_0$, and therefore in the box $B_0$ about $y_0$ with side $6\eps$.  Since they are grid points in the box $B_0$, they can be connected by a grid path in $B_0$, which by assumption is in $\omhp$.

Now we have a grid path from $x_0$ to $\xtw_1$, with
$|\xtw_1 - x(\tau)| < \eps$.  We proceed by induction, increasing $t$ by $\tau$ at each step:  Assume for $j \geq 1$ we have a path of grid points from $x_0$ to $\xtw_j$ with $|\xtw_j - x(j\tau)| < \eps$.
We choose $y_j \in Y$ so that
$|y_j - x(j\tau)| < \eps$.  Then $|\xtw_j - y_j| < 2\eps$, and as before
$|x((j+1)\tau) - y_j| < |x((j+1)\tau) - x(j\tau)| + |x(j\tau) - y_j| <
2\eps$.
 We choose a grid point $\xtw_{j+1}$ with $|\xtw_{j+1} -x((j+1)\tau)| < \eps$, so that $|\xtw_{j+1} - y_j| < 3\eps$.  Then 
$\xtw_j$ and $\xtw_{j+1}$ are both in the box $B_j$ about $y_j$, and they
can be connected by a grid path in $B_j \subseteq \omp$, completing the induction step.  Making
about $1/\tau$ steps, we connect to $x(1) = z$ at the last step.

Finally if $x \in \omhp$ is within distance $\rho$ of $\Gamma$, it is in one of the sets
${\cal N}_z$ previously chosen and has the form $x = (x^1,\gamma(x^1) - y)$ with $y \leq C_0\rho$,
using (4.1).  Then $b = (x^1,\gamma(x^1) - 2C_0\rho) \in {\cal N}_z$,
and $\mbox{dist}(b,\Gamma) \geq 2\rho$, again by (4.1).  Provided $h < C_0\rho$, we
can connect $x$ by a vertical grid path to a point in $\om_\rho\cap\omhp$.  In this way
all points in $\omhp$ can be connected.

\bigskip

{\bf Proof of Lemma 2.3.}
We prove (2) in $\R^2$ with remarks about $\R^3$ and then comment on (3). 
We note that if $u$ has mean value zero and
$c$ is any constant, then $\sum (u + c)^2\,h^2  \geq \sum u^2\,h^2$, and we may add a constant
to $u$ as needed in proving the inequality.  We prove it first on a subdomain away from $\Gamma$,
using the standard Poincar\'{e} inequality with mean value zero, and then take care of points near $\Gamma$ by summing differences.

With $\rho > 0$ and the sets ${\cal N}_z$ as chosen previously, let
$\om_\rho = \{x \in \omp: \mbox{dist}(x,\Gamma) > \rho\}$
and similarly for $\rho'' < \rho' < \rho$.  
Let $\zeta$ be a cut-off function such that $\zeta = 1$ within distance $\rho$ of $\Gamma$ and supported within distance $2\rho$. 
We will prove the inequality on $\om_\rho$ and then for $\zeta u$ near $\Gamma$.
We assume $h$ is small relative to $\rho$.

We will interpolate the grid function $u$ to $\utw$ on $\om_{\rho'}$ using bilinear (or trilinear) interpolation on grid squares (or cubes) in $\R^2$ (or $\R^3$).
With $d = 2$, we will first show that
\beq \int_{\om_{\rho'}} |\nabla \utw|^2 \,dx \leq  
        C \sum_{x \in \om_{\rho''}}  |D_\nu^2 u(x)|^2 \,h^2 \eeq
and
\beq \sum_{x \in \om_{\rho}}  u(x)^2 \,h^2 
                 \leq C \int_{\om_{\rho'}} \utw^2 \,dx \eeq
To verify these, suppose $x_{0,0}, x_{1,0}, x_{0,1}, x_{1,1}$ are vertices of a grid square in
$\om_{\rho''}$.  The bilinear interpolation is
$$ \utw(x^1,x^2) \eq a_{00} +  a_{10}x^1 +  a_{01}x^2 + a_{11}x^1x^2/h  $$
where $a_{00} = u_{0,0}$ and, for $(\alpha,\beta) \neq (0,0)$, $a_{\alpha\beta}$ is a sum of the differences $Du$ on the square.  If $Q = \max{|Du|}$ for the four differences, then
$|\nabla \utw| \leq CQ$, leading to an estimate for the integral on the square
from which (4.2) follows by summing over squares.
For the second inequality, suppose $U = |u_{0,0}|$ is the largest of the four $|u_{i,j}|$
on a square.
For each $ (\alpha,\beta)$ we have $|a_{\alpha,\beta}| \leq CU/h$ and 
$|\utw(x)| \geq U(1 - C_1|x|/h)$.  Then for $\theta$ small enough, depending on $C_1$,
$$ \int_0^{\theta h} \int_0^{\theta h} \utw^2 \,dx^1dx^2 \geq 
      U^2h^2 \int_0^{\theta} \int_0^{\theta}(1 - C_1|s|)^2\,ds^1ds^2 
          \geq c_0  U^2h^2  $$
and (4.3) follows.  A very similar argument works in $\R^3$. 

Now suppose we adjust $\utw$ by a constant so that it has mean value zero on $\om_\rho'$
and change $u$ by the same constant.  A standard form of the Poincar\'{e} inequality gives
\beq \int_{\om_{\rho'}} \utw^2 \,dx 
           \leq C  \int_{\om_{\rho'}} |\nabla \utw|^2 \,dx  \eeq
Then, combining the three inequalities, we have the desired conclusion for the subdomain,
\beq \sum\nolimits_{x \in \om_{\rho}}  u(x)^2 \,h^2 
                 \leq C \sum\nolimits_{x \in \om_{\rho''}}  |D_\nu^2 u(x)|^2 \,h^2 \eeq

Finally we prove the inequality for $\zeta u$ in the set ${\cal N}_z$.  For
$x^1$ such that $(x^1,\gamma(x^1)) \in {\cal N}_z$, let
$b(x^1) =  \gamma(x^1) - 2C_0\rho$.  By (4.1), $(x^1,b(x^1))$ is at distance at least
$2\rho$ from $\Gamma$, so that $\zeta(x^1,b(x^1)) = 0$.  If $(x^1,x^2) \in {\cal N}_z\cap\omhp$
is within $\rho$ of $\Gamma$, then $x^2 \geq \gamma(x^1) - C_0\rho$, again by (4.1).
Now for all grid points $(x^1,x^2)$ with $b(x^1) \leq x^2 \leq \gamma(x^1)$
we write $(\zeta u)(x^1,x^2)$ as a sum of vertical differences and use the Cauchy-Schwarz 
inequality to get
$$  |(\zeta u)(x_1,x_2)|^2 \leq
                    C \sum\nolimits_{s \leq \gamma(x_1)} |D_2(\zeta u)(x_1,s)|^2\,h  $$
We sum over $x_2$ and then $x_1$ to obtain
\beq \sum\nolimits_{\omhp\cap{\cal N}_z} |(\zeta u)(x_1,x_2)|^2 \,h^2   
                \leq C \sum\nolimits_{\omhp\cap{\cal N}_z} |D_2(\zeta u)(x_1,x_2)|^2\,h^2  \eeq
Since each point in $\omhp$ within $\rho$ of $\Gamma$ is in some ${\cal N}_z$ and $\zeta = 1$ for
such points, the estimate (2.7) is complete.

The proof of (3) is similar, using the fact that for $\utw$ on $\bbox - \omp$ with 
$\utw = 0$ on $\pa\bbox$, the inequality like (4.4) holds on the outer domain.
This is one standard form of the Poincar\'{e} inequality; it can easily be proved
by contradiction using the Rellich Compactness Theorem.

\bigskip

{\bf Proof of Lemma 2.6.}
We prove the first statement in $\R^2$ and remark on the other cases.
We need to use the assumption that the cut point for each interval in $\intshcut$ is at least
distance $\delta h$ away from either endpoint, with $\delta > 0$ fixed.
This implies $\delta \leq \xi^\nu \leq 1 - \delta$, and it is equivalent to show that
\beq \sum\nolimits_{\intshp \cup \intshcut} (D_\nu v^+)^2\,h^2 \leq C_1  
      \sum\nolimits_{\intshm \cup \intshcut} (D_\nu v^-)^2\,h^2 \equiv E \eeq
We will actually show that there is some $w \in \fcnsp$, not assumed harmonic, with
$Mw = \varphi = Mv^-$ so that the inequality holds with $w$ in place of $v^+$.  Then
according to Lemma 2.5, $\langle v^+,v^+ \rangle^+ \leq \langle w,w \rangle^+$, and the
result follows.   

From the discrete Poincar\'{e} inequality, Lemma 2.3(3),
we have $\sum (v^-)^2\,h^2 \leq CE$ on $\omhm$.  We use a partition of
unity to write $v^-$ as a sum of terms in $\fcnsm$ each of which is supported in a set where
$\Gamma$ has the form $x^2 = \gamma(x^1)$, or the reverse, where $x = (x^1,x^2)$. With a
smooth cut-off function $\zeta$, and for an interval $[x_{i,j},x_{\ip,j}]$ with
$x_{i,j} \in \omhm$, we have
$D(\zeta v^-)_{\iph,j} = \zeta_{\ip,j}(Dv^-)_{\iph,j} + (D\zeta)_{\iph,j}v^-_{i,j}$,
from which we see that 
$$  \sum\nolimits_{\intshm \cup \intshcut} (D_\nu (\zeta v^-))^2\,h^2 \leq CE  $$
It is enough to find $w$ for such a localized function.  Thus we replace $v^-$ with
a function $v$ supported in a set where $x^2 = \gamma(x^1)$ on 
$\Gamma$, $x^2 < \gamma(x^1)$ in $\omm$, and $x^2 > \gamma(x^1)$ in $\omp$. 

We will define $w$ first by even reflection across $\Gamma$ along vertical lines and then find extended values in $\omhm$ at endpoints of horizontal intervals.  It is helpful to note that
an extended value affects only the interval where it is defined.
In the first part $w$ at $x = (x^1,x^2)$ will be
determined by $v$ near $R(x) =(x^1,2x^* - x^2)$, where $x^* = x^*(x^1)$ and $(x^1,x^*) \in
\Gamma_h^1$.  We construct $w$ in such a way that for each $I \in \intshp \cup \intshcut$,
$|Dw| \leq C \sum (|Dv| + |v|)$, where $|Dw|$ is the difference on $I$.  Here the sum is over intervals within distance $c_0h$ of $R(x)$, where $x$ is one endpoint of $I$, and  $c_0$ is some constant depending on a bound for $|\gamma'|$.
A given interval in $\intshm \cup \intshcut$ will occur in such a sum only for a bounded number of $I$.  The needed estimate
for $w$ will then follow, and the proof of (4.7) will be complete.  

With $x^1$ fixed, we define $w$ as a function of $x^2$.  For
convenience we assume $0 \leq \gamma(x^1) \leq h$ and $(x^1,x^*) \in \Gamma_h^1$,
with $x^* = (1-\theta)h$ and $\delta \leq \theta \leq 1-\delta$.
We temporarily write $v$ and $w$ as functions of $x \in \R$ rather than $x^2$.  We can use
piecewise linear interpolation to extend $v$ from grid points to all $x \leq h$.
We will write $v_j = v(jh)$, $j\leq 1$.  
For $j \geq 1$ we will define $w_j = w(jh)$ by even
reflection of $v$ about $x^*$, that is
\beq w_j \eq v(2x^* - jh) \eq v((2 - 2\theta -j)h)\,, \qquad j \geq 1 \eeq
while for $j = 0$ we need a special definition so that $Mw = \varphi$.  We need to express
differences of $w_j$ in terms of those for $v_j$ to verify the boundedness condition.  We separate into two cases:  (1) $\theta \geq \lilhalf$ and  (2) $\theta \leq \lilhalf$.
For case (1), $0 \leq 2-2\theta \leq 1$ and
 \beq w_j \eq (2\theta -1)v_{-j} \p (2-2\theta) v_{-j+1}  \,, \qquad j \geq 1 \eeq
Then
\beq w_{\jp} - w_j \eq (2\theta -1)(v_{-j-1} - v_{-j}) + (2-2\theta) (v_{-j} - v_{-j+1}) \eeq
so that $w_{\jp} - w_j$ is an interpolation of differences of $v$ within $\omhm$, $j\geq 1$.
We define $w_0$ from the requirement
$Mw = g = Mv$, or  $\theta w_0 + (1-\theta)w_1 = \theta v_0 + (1-\theta)v_1$ and
\beq w_0 - w_1 \eq (v_0 - w_1) + \theta^{-1}(1-\theta)(v_1 - w_1)  \eeq
In case (1) the two differences on the right are combinations of differences
appearing in $E$ since $w_1$ is an interpolation of $v_0$ and $v_{-1}$. 

For case (2), $\theta \leq \lilhalf$, $1 \leq 2-2\theta \leq 2$, so that
 \beq w_j \eq 2\theta v_{-j+1} \p (1-2\theta) v_{-j+2}  \,, \qquad j \geq 1 \eeq
and
\beq w_{j+1} - w_j \eq 2\theta(v_{-j} - v_{-j+1}) + (1-2\theta) (v_{-j+1} - v_{-j+2}) \eeq
Defining $w_0$ again so that $Mw = Mv$ and using (4.11),(4.12) we find that 
$w_0 - w_1 = v_1 - v_0$. 

We have now defined $w$ at all grid points in $\om^+$ and at endpoints in $\om^-$ of vertical cut intervals.  We denote the latter values as $w^{(2)}$ to indicate the interval.
We need to assign values $w^{(1)}$
at the endpoints in $\om^-$ on horizontal intervals satisfying the requirement $Mw = \varphi$.
It may be that one endpoint in $\om^-$ is assigned multiple values $w^{(1)}, w^{(2)}$.

We consider a horizontal interval $[x_{i,j},x_{\ip,j}] \in \intshcut$ with 
$x_{i,j} \in \omhm$, $x_{\ip,j} \in \omhp$, and 
$\theta_1x_{i,j} + (1-\theta_1)x_{\ip,j} \in \Gamma_h^1$, where $\delta \leq \theta_1 \leq 1$
by assumption.
We have already assigned $w_{\ip,j}$ and we need to define $w_{i,j}^{(1)}$.
From $Mw = \varphi = Mv$  we find that $w_{i,j}^{(1)}$ is defined by
\beq w_{i,j}^{(1)} - w_{\ip,j} \eq (v_{i,j} - w_{\ip,j})
                + \theta_1^{-1}(1-\theta_1)(v_{\ip,j}^{(1)} - w_{\ip,j}) \eeq
We need to see that the two terms on the right are appropriately bounded.
Note that $w_{\ip,j}$ was defined by interpolation from $v_{\ip,k}$ with $k -j$ bounded, depending on $|\gamma'|$.  We can form a path of grid intervals
$[x_{i,j},x_{i,\ell},x_{\ip,\ell},x_{\ip,k}]$ with $\ell < j$, $\ell \leq k$ and $\ell - j$ bounded in terms of $|\gamma'|$,
such that each endpoint is in $\omhm$,
except possibly $x_{\ip,k} = x_{\ip,j}$.  By adding and subtracting, $v_{i,j} - w_{\ip,j}$ is
then a sum of differences of $v$ occurring in $E$.  For the second term, $\theta_1^{-1} \leq \delta^{-1}$, and
we can use the same path with the additional term $v_{\ip,j}^{(1)} - v_{i,j}$.  The reverse case with $x_{i,j} \in \omhp$, $x_{\ip,j} \in \omhm$ is similar.
It is possible that $x_{i,j}$ is an endpoint of two horizontal cut intervals and
is assigned values by each.

Finally we need to check that the difference of $w$ on a horizontal interval
$[x_{i,j},x_{\ip,j}]$ with both endpoints in $\omhp$ can be written as differences of $v$
near a reflected point. 
The value of $w$ at each point is a vertical interpolation of values of $v$ at grid points in
$\omhm$, and possibly an extended value in $\omhp$. If $k$ is an integer near
$R(x_{i,j}^1)$, we can write $w_{i,j}$ as
$v_{i,k}$ plus vertical differences near $k$.  Similarly $w_{\ip,j}$ equals $v_{\ip,\ell}$ plus differences, with $\ell - k$ bounded. We can write $v_{\ip,\ell} - v_{i,k}$ as a sum of differences on intervals as we did above.

For domains in $\R^3$ we can use the same argument with two horizontal directions.  To estimate
$v^-$ by $v^+$, we use Lemma 2.3(2).

\ssection{Numerical examples}

We present three examples of the method described in Sec. 1.  In each case, with the quadratic interpolation (1.6), we find the solution of the Dirichlet problem has accuracy $O(h^2)$.  Furthermore the accuracy near the 
boundary $\Gamma$ is nearly $O(h^3)$, a property that is characteristic of the
Shortley-Weller method \cite{matyam,weynans,YMSW}.  We also test the version with the linear interpolation (1.9).
We find it is $O(h^2)$ accurate but less accurate than the Shortley-Weller method.

Given a curve $\Gamma$, we choose an exact harmonic function $u$ and
specify the boundary value $g$ on $\Gamma$.  We compute the solution
of the Dirichlet problem in Matlab
according to the method of Sec. 1.  We use the computational box
$\B = (-1,1)\times (-1,1)$ and choose $h = 2/N$ with integer $N$.
We identify the cut intervals and intersection points.  We solve the
equation $A^h\varphi^h = g^h$ using GMRES.  For each provisional $\varphi^h$
we solve (1.4), (1.5) for $u^h$ by inverting $\laph$, and then interpolate to find
$f^h = A^h\varphi^h$.  After enough iterations so that the residual
$g^h - A^h\varphi^h$ is
within a specified tolerance, we compare the computed harmonic function
$u^h$ in $\omhp$ with the exact solution $u$.

{\bf Example 1.}
We choose $\Gamma$ to be the circle $x_1^2 + x_2^2 = r^2$ with $r = .6$, and
$u(x_1,x_2) = \sin(x_1 + x_2)\exp(x_1 - x_2)$.  We solve with either linear or quadratic interpolation.   Results are shown in Table 5.1.  In each row we
display $N$; the specified error tolerance for the relative residual, such as
e$-8 = 10^{-8}$; the number of iterations required in GMRES;
the maximum error in $u^h$ on the interior domain $\omhp$; the maximum error within distance $4h$ of the boundary; and the same for distance $2h$.  For the linear version the maximum near the boundary is the same as the maximum in the domain.
With linear interpolation we see convergence close to $O(h^2)$, as
expected.
Using quadratic interpolation, the maximum error is more clearly $O(h^2)$ and the errors are much smaller than in the linear case.  The errors near the boundary have higher accuracy, and we see about $O(h^3)$ accuracy for the
largest $N$, provided we reduce the tolerance.  This is strong evidence that the solution obtained is the Shortley-Weller solution.  The number of iterations needed in GMRES
increases with $N$.  We have chosen the tolerance so that reducing it further does
not change the errors significantly.

\begin{table}[ht]
\caption{Errors for a circle with radius .6} 
\vspace{6pt}
\centering 
\begin{tabular}{| c | c | c | c | c | c | c |} \hline
method & N & tol & it'ns & err$(\omhp)$ & err$(4h)$ & err$(2h)$ \\ \hline  \hline
\multirow{4}{*}{linear} 
& 64 & e-6 &  13 & 3.10e-4 & same & same \\ \cline{2-7}
& 128 & e-6 & 16 & 8.60e-5 &" & "\\ \cline{2-7}
& 256 & e-6 & 19 & 2.30e-5 &" & "\\ \cline{2-7}
& 512 & e-6 & 23 & 5.76e-6 & "& "\\ \hline
\hline
\multirow{7}{*}{quadratic}
& 32 & e-8 & 17	& 4.93e-5	& 4.93e-5&     4.13e-5  \\ \cline{2-7}
& 64 & e-8 & 23	& 1.06e-5	& 8.56e-6&     5.39e-6  \\ \cline{2-7}	
& 128 & e-8 & 26 & 2.59e-6	& 1.28e-6 &     7.58e-7  \\ \cline{2-7}
& 256 & e-8 & 31 & 6.34e-7   & 1.76e-7& 	1.04e-7  \\ \cline{2-7}
& 512 & e-8 & 46 & 1.56e-7   & 2.82e-8&      2.82e-8  \\ \cline{2-7}
& 256 & e-9 & 38 & 6.34e-7   &1.76e-7 &	1.04e-7 \\ \cline{2-7}
& 512 & e-9 & 52 & 1.56e-7   &2.25e-8 &	1.29e-8    \\ \hline
\end{tabular}
\end{table}

{\bf Example 2.}
For our second example $\Gamma$ is the ellipse $x_1^2/a^2 + x_2^2/b^2$ with 
$a = .6$, $b = .4$, and the exact solution $u$ is the same as before.
Table 5.2 is like Table 5.1,
except that in the last column we display the maximum error at interior irregular points, i.e., grid points which are the endpoints of cut intervals.  The results are similar
to those for the circle, with second order accuracy for both methods, smaller
errors with quadratic interpolation, and about $O(h^3)$ accuracy at the irregular points.

\begin{table}[ht]

\caption{Errors for an ellipse} 
\vspace{6pt}
\centering 
\begin{tabular}{| c | c | c | c | c | c |} \hline
method & N & tol & it'ns & err$(\omhp)$ & irreg err \\ \hline  \hline
\multirow{4}{*}{linear} 

& 64 & e-6 &  15 & 3.35e-4 & same  \\ \cline{2-6}
& 128 & e-6 & 20 & 8.50e-5 &" \\ \cline{2-6}
& 256 & e-6 & 34 & 2.14e-5 &" \\ \cline{2-6}
& 512 & e-6 & 34 & 5.42e-6 & "\\ \hline
\hline
\multirow{6}{*}{quadratic}
& 32 & e-8 & 23	&  3.90e-5 &	3.26e-5 \\ \cline{2-6}
& 64 & e-8 & 26	&  8.00e-6 & 4.51e-6\\ \cline{2-6}	
& 128 & e-8 & 32 &  1.90e-6	 & 5.83e-7\\ \cline{2-6}
& 256 & e-8 & 57 &  4.62e-7	& 7.56e-8\\ \cline{2-6}
& 512 & e-8 & 70 &   1.14e-7 &  2.06e-8\\ \cline{2-6}
& 512 & e-9 & 80 &   1.14e-7	& 9.67e-9\\ \hline
\end{tabular}
\end{table}

{\bf Example 3.}
For our final example $\Gamma$ is the Cassini oval
$$ \left( (x_1-a)^2 + x_2^2\right)\left( (x_1+a)^2 + x_2^2\right) = b^4 $$
with $a = .4$ and $b = .44$, which is partly concave.  The harmonic function
is $u = \log r$, where $r^2 = (x_1-1)^2 + (x_2-1)^2$.  The results, shown in Table 5.3, are similar to the previous ones, though the orders of accuracy are not quite as close.

\begin{table}[ht]
\caption{Errors for the Cassini oval} 
\vspace{6pt}
\centering 
\begin{tabular}{| c | c | c | c | c | c |} \hline
method & N & tol & it'ns & err$(\omhp)$ & irreg err\\ \hline  \hline
\multirow{6}{*}{linear} 
& 32 & e-6 &  14 & 1.94e-4 & same  \\ \cline{2-6}
& 64 & e-6 & 19 & 3.22e-5 &" \\ \cline{2-6}
& 128 & e-6 & 30 & 1.23e-5 &" \\ \cline{2-6}
& 256 & e-6 & 27 & 3.71e-6 & "\\ \cline{2-6}
& 512 & e-6 & 32 & 1.99e-6 & "\\ \cline{2-6}
& 512 & e-8 & 62 & 1.01e-6 & "\\  \hline
\hline
\multirow{7}{*}{quadratic}
& 32 & e-8 & 21	&  2.18e-5 &	1.37e-5 \\ \cline{2-6}
& 64 & e-8 & 35	&  6.20e-6 & 1.50e-6\\ \cline{2-6}	
& 128 & e-8 & 44 &  1.64e-6	 & 2.47e-7\\ \cline{2-6}
& 256 & e-8 & 52 &  4.23e-7	& 3.15e-8\\ \cline{2-6}
& 512 & e-8 & 80 &   1.07e-7 &  9.53e-9\\ \cline{2-6}
& 256 & e-9 & 58 &   4.23e-7 &  3.32e-8\\ \cline{2-6}
& 512 & e-9 & 87 &   1.07e-7	& 4.39e-9\\ \hline
\end{tabular}
\end{table}

The Shortley-Weller method is awkward to implement directly since it requires a
matrix particular to the domain and the matrix is not symmetric.  The indirect procedure here has the advantage that each iteration requires only the inversion of
$\laph$ on a rectangular grid.  However we have seen that the number of iterations needed grows with the refinement.  This disadvantage may be related to the difference between the norms $\|\cdot\|_1$ and $\|\cdot\|_2$ described in the Main Theorem
and in Sec. 3.  In contrast, the method of \cite{YWjcp} uses an equilibration procedure for function values at the intersection points, and the number of iterations is almost independent of refinement.  However we have not been able to extend the analysis presented here to include this equilibration.

\bibliographystyle{amsplain}

\end{document}